\documentclass[10pt,a4paper]{article}
\usepackage[latin2]{inputenc}
\usepackage[english]{babel}
\usepackage{a4wide}
\usepackage[english]{babel}
\usepackage[nottoc]{tocbibind}
\usepackage{amsmath,amsfonts,amstext,amsbsy,amsopn,amssymb,amsthm,amscd}
\usepackage{amsxtra,latexsym}
\usepackage{bbm}
\usepackage{exscale}
\usepackage{paralist}
\usepackage{color,graphics}
\usepackage[article]{math_env} 
\numberwithin{equation}{section}

\newcommand{\field}[1]{\mathbb{#1}}
\newcommand{\N}{\field{N}}
\newcommand{\R}{\field{R}}

\def\bbA{{\mathbb A}}  \def\bbC{{\mathbb C}}

  \def\bbL{{\mathbb L}}

\def\calP{{\mathcal P}} \def\calQ{{\mathcal Q}} 
\def\calS{{\mathcal S}} \def\calT{{\mathcal T}}

 \newcommand{\dx}{\,\mathrm{d} x}
 \newcommand{\dy}{\,\mathrm{d} y}

\newcommand{\abs}[1]{\left\lvert#1\right\rvert}      
\newcommand{\norm}[1]{\left\lVert#1\right\rVert}

\newcommand{\wt}{\widetilde}
\DeclareMathOperator{\grad}{grad}
\DeclareMathOperator{\supp}{supp}
\DeclareMathOperator{\Lin}{Lin}

\DeclareMathOperator{\curl}{curl}
\DeclareMathOperator{\sym}{sym}
\renewcommand{\nabla}{\mathrm{D}}
\newcommand{\nquad}{\mspace{-18.0mu}}
\DeclareMathOperator{\Curl}{Curl}
\DeclareMathOperator{\dyw}{div}
\DeclareMathOperator{\Dyw}{Div}
\DeclareMathOperator{\Sym}{Sym}

\newcommand{\D}{{\mathbb C_{{\rm e}}}}

\newcommand{\Ha}{{\mathbb C_{{\rm micro}}}}
\newcommand{\Lc}{{\mathbb L}_{{\rm c}}}

\newcommand{\idd}{\mathbbm{1}}

\begin{document}
\title{A local regularity result for the relaxed micromorphic model based on inner variations}
\author{%
Dorothee Knees
\thanks{
Institute of Mathematics, University of Kassel, Heinrich-Plett Str. 40, 34132 Kassel, Germany, dknees@mathematik.uni-kassel.de
}
\and
Sebastian Owczarek\thanks{Faculty of Mathematics and Information Science, Warsaw University of Technology, ul. Koszykowa 75, 00-662 Warsaw, Poland, s.Owczarek@mini.pw.edu.pl}
\and
Patrizio Neff\thanks{Lehrstuhl f\"ur Nichtlineare Analysis und Modellierung, Fakult\"at f\"ur Mathematik, Universit\"at Duisburg-Essen, Campus Essen, Thea-Leymann Str. 9, 45127 Essen, Germany, patrizio.neff@uni-due.de}
}

\date{ }

\maketitle

\begin{abstract}
In this paper we study local higher regularity properties of a linear elliptic system that is coupled with a system of Maxwell-type. The regularity result is proved by means of a modified finite difference argument. 
These  modified finite differences are based on inner variations combined with a Piola-type transformation in order to preserve the $\curl$-structure in the Maxwell system. The result is applied to the relaxed micromorphic model. 
\end{abstract}




%
\section{Introduction}
\label{s:introduction}

In this paper we study local higher regularity properties of a linear elliptic system that is coupled with a system of Maxwell-type. As an application we have in mind the linear relaxed micromorphic model developed in  \cite{NGperspective,Ghibastatic,Ghibadyn,Owczghibaneffexist,exprolingalberdi}. The latter is formulated on a bounded domain $\Omega\subset \R^3$ and reads as follows: given positive definite and symmetric material dependent coefficient tensors $\bbC_e\in \Lin(\Sym(3),\Sym(3))$, $\bbC_\text{micro}\in \Lin(\Sym(3),\Sym(3))$ and $\bbL_c\in \Lin(\R^{3\times 3},\R^{3\times 3})$ determine a displacement field $u:\Omega\to\R^3$ and a non-symmetric micro-distortion tensor $P:\Omega\to\R^{3\times 3}$ satisfying 
\begin{align}
\label{intro:e1}
 0&=\Dyw\Big(\bbC_e\sym(\nabla u - P)\Big) + f\quad \text{in }\Omega,\\
 0&= -\Curl\Big(\bbL_c\Curl P\big) + \bbC_e\sym(\nabla u - P) - \bbC_\text{micro}\sym P + M\quad \text{in }\Omega\,.
 \label{intro:e2}
\end{align}
Here, $f:\Omega\to\R^3$ is a given volume force density and $M:\Omega\to \R^{3\times 3}$ a body moment tensor. These equations need to be closed by boundary conditions. We refer to Section \ref{sec:micro} for more information about the model and to Section \ref{s:notation} for details concerning the notation. 

Weak solutions to this model satisfy $u\in H^1(\Omega)$ and $P\in H(\Curl;\Omega)=\{Q\in L^2(\Omega)\, |\,\,\Curl Q\in L^2(\Omega)\}$. Higher regularity properties of linear systems can be obtained by many different arguments. We present here a local regularity result that is based on a modified finite difference argument. Assuming enough regularity for the data $f$ and $M$ we show the higher local regularity  $u\in H^2_\text{loc}(\Omega)$ and $P,\Curl P\in H^1_\text{loc}(\Omega)$, see Theorem \ref{them_reg_cube}.  

The obtained regularity result for the tensor $P$, in particular $P\in H^1_\text{loc}(\Omega)$, is at the same time surprising and not surprising: The bilinear form associated with \eqref{intro:e1}--\eqref{intro:e2} reads 
\begin{multline}
\label{intro:bili_micromorph}
a((u,P);(v,Q)):=\int_\Omega \langle\bbC_e\sym(\nabla u- P),\sym(\nabla v - Q)\rangle +
\langle \bbC_\text{micro}\sym
P,\sym Q \rangle
\\+\langle \bbL_c\Curl P, \Curl Q\rangle\,\dx\,.
\end{multline}
From the positivity of the involved material tensors one then finds 
\begin{align}
 a((u,P);(u,P))\geq c\Big(\norm{\sym\nabla u}^2_{L^2(\Omega)} +\norm{\sym P}_{L^2(\Omega)}^2 
 +\norm{\Curl P}_{L^2(\Omega)}^2\Big)\,.
\end{align}
As a surprising fact it is shown in \cite{neff2012poincare,NeffLP,lewintan_neff_2021,LMN21} that on suitable subspaces of $H^1(\Omega)\times H(\Curl;\Omega)$ the right hand side can further be estimated by generalized versions of Korn's inequality as follows:
\begin{align}
\label{intro:Korn}
 a((u,P);(u,P))\geq c_K\Big(\norm{ u}_{H^1(\Omega)}^2 + \norm{P}_{L^2(\Omega)}^2 + \norm{\Curl P}^2_{L^2(\Omega)}\Big)\,.
\end{align}
Hence, in particular the $L^2$-norm of the tensor $P$ and not only of its symmetric part is controlled. This coercivity estimate (see also assumption (A1) here below) is the starting point for the regularity arguments. By choosing well-adapted test functions based on modified finite differences one ultimately obtains an additional derivative (locally) in each of the terms on the right hand side of \eqref{intro:Korn}.

Let us give more details about the arguments developed in this paper: 
Standard proofs based on finite differences to obtain higher local regularity for elliptic PDEs (in the sense of higher differentiability in the scale of standard Sobolev spaces) start from the weak formulation of the PDE and work with test functions of the type $u_h(x):= \triangle^h(\varphi^2\triangle_h u)$, where $u\in H^1(\Omega)$ is a weak solution, $\varphi\in C_0^\infty(\Omega)$ is a cut-off function, $h\in \R^3$ and $\triangle_h v(x):= v(x+h)-v(x)$, $\triangle^h v(x)= v(x)-v(x-h)$. If $\abs{h}>0$ is small enough then $u_h\in H_0^1(\Omega)$, and hence, $u_h$ is an admissible test function for the weak formulation. 
With the help of this test function one next shows that for every $\wt\Omega\Subset\Omega$ there is a constant $c>0$ such that $\abs{h}^{-1}\norm{u(\cdot + h) - u(\cdot)}_{H^1(\wt\Omega)}\leq c $  for small enough $h\in \R^3$. From this, by \cite[Lemma 7.24]{GiTr86} 
one finds $u\in H^2_\text{loc}(\Omega)$. To obtain such an estimate, the ellipticity/coercivity of the associated bilinear form is exploited and terms of lower order are estimated from above yielding suitable powers of $\abs{h}$. Transferring this procedure to a bilinear form involving the $\Curl$-term we are faced with the following problem that we explain here in a simplified setting: For vectors $p,q\in H_0(\curl;\Omega,\R^3)$ (the completion of $C_0^\infty(\Omega,\R^3)$ with respect to the norm $\norm{p}_{H(\curl,\Omega)}=\norm{p}_{L^2(\Omega)} + \norm{\curl p}_{L^2(\Omega)}$)   let 
\begin{align*}
 a(p;q):=\int_\Omega \langle\curl p, \curl q\rangle\,\dx\, 
\end{align*}
and assume that $p\in H_0(\curl;\Omega)$ satisfies $a(p,q)=\int_\Omega \langle M ,q \rangle\,\dx$ for some smooth right hand side $M:\Omega\to \R^3$ and all $q\in H_0(\curl;\Omega)$. 
With $q_h=\triangle^h(\varphi^2\triangle_h p)$ as above we find: $\curl q_h = \triangle^h(\varphi^2\triangle_h\curl p + (\grad \varphi^2)\times \triangle_h p)\in L^2(\Omega)$ and 
\begin{align*}
 a(p;q_h) &= \int_\Omega \langle\curl p, \triangle^h(\varphi^2\triangle_h\curl p + (\grad \varphi^2)\times \triangle_h p)\rangle\,\dx\\
 &= \int_\Omega \varphi^2 \abs{\triangle_h \curl p}^2\,\dx 
 + 2 \int_\Omega  \big\langle\varphi \triangle_h \curl p,\big(\grad \varphi\times \triangle_h p\big)\big\rangle\,\dx \,.
\end{align*}
Assume now that we have $\abs{\int_\Omega \langle M, q_h\rangle,\dx }\leq c\abs{h}^2\norm{p}_{H(\curl;\Omega)}$ with some $M$-dependent constant $c>0$, and hence,  
\begin{align}
\label{eq:introe2}
\int_\Omega \varphi^2 \abs{\triangle_h \curl p}^2\,\dx 
 + 2 \int_\Omega  \big\langle\varphi \triangle_h \curl p,\big(\grad \varphi\times \triangle_h p\big)\big\rangle\,\dx\leq c\abs{h}^2\norm{p}_{H(\curl;\Omega)}\,. 
\end{align}
An estimate of the type 
\[
 \int_\Omega \varphi^2 \abs{\triangle_h \curl p}^2\dx \leq c \abs{h}^2
\]
implies that  $\curl p\in H^1_\text{loc}(\Omega)$. In order to obtain such an estimate, by Young's inequality, the first factor in the second term on the left hand side of \eqref{eq:introe2} can be absorbed in the first term so that we ultimately  obtain
\begin{align}
 \int_\Omega \varphi^2 \abs{\triangle_h \curl p}^2\dx \leq c\,\big(\norm{\grad \varphi\times \triangle_hp}_{L^2(\Omega)}^2 + \abs{h}^2\norm{p}_{H(\curl;\Omega)}\big)\,.
\end{align}
The problem now is that the term $\norm{\grad \varphi\times \triangle_hp}_{L^2(\Omega)}$ can neither be absorbed on the left hand side nor be estimated by $c\abs{h}$ since the latter would require $p\in H^1(\Omega)$, which we do not have. 

The main idea  to circumvent this problem is to replace standard finite differences with finite differences that are 
based on inner variations and that involve a modified Piola transformation for $p$, see formula  \eqref{def:Piolatypetrafo} and Lemma \ref{lem:piolacurl}.  
For such test functions, terms of the type $\norm{\grad\varphi\times \triangle_hp}_{L^2(\Omega)}$ do not appear. 

In the context of regularity theory for nonlinear (and nonsmooth) elliptic systems, inner variations  proved to be very powerful already in many different contexts and we mention here \cite{Nes06,NesAlberreg} for regularity results in visco-plasticity, 
\cite{Neffkneesreg} for a global regularity result for weak solutions to a quasilinear elliptic equation emerging in a time-discretized plasticity-model, and \cite{KS12} for the global regularity of minimizers of uniformly convex but not necessarily smooth functionals involving for example friction and being defined on domains with cracks and contact conditions. 
As already mentioned, for linear systems there are also many other ways to prove higher regularity of solutions.  The arguments shown in this paper develop their full power in nonlinear problems and for problems formulated on domains with Lipschitz boundaries,  \cite{KNO22-glob}.  

The paper is organized as follows: In Section \ref{s:notation} we clarify the notation and introduce the spaces under consideration. Section \ref{s:regularity} is devoted to the formulation of the main result, Theorem \ref{them_reg_cube}, and its proof via the modified finite differences. Finally, in Section 4 we apply the result to the above mentioned linear relaxed micromorphic model and discuss consequences for FEM approximations.

\section{Notation}  
\label{s:notation}

For  vectors $a,\,b\in\R^n$, we define the scalar product $\langle a,b\rangle:=\sum_{i=1}^n a_ib_i$, the Euclidean norm $\abs{a}^2:=\langle a,a\rangle$ and the dyadic product $a\otimes  b=(a_ib_j)_{i,j=1}^n\in\R^{n\times n}$, where $\R^{m\times n}$ will denote the set of real $m\times n$ matrices. The identity tensor on $\R^{n\times n}$ will be denoted by $\idd$. For  matrices $P,Q\in \R^{m\times n}$, we 
define the standard Euclidean scalar product $\langle P,Q\rangle:=\sum_{i=1}^m\sum_{j=1}^n P_{ij}Q_{ij}$ and the Frobenius-norm $\|P\|^2:=\langle P,P\rangle$. $P^T\in \R^{n\times m}$ denotes the transposition of the matrix $P\in \R^{m\times n}$ and for $P\in \R^{n\times n}$, the symmetric part of $P$
will be denoted by $\sym P=\frac{1}{2}(P+P^T)\in \Sym(n)$. We recall the following scalar product rule for matrices $P\in \R^{n \times k}$, $Q\in \R^{k \times m}$ and $R\in \R^{n \times m}$: 
\begin{align}
\label{eq:sp_matrix}
\langle PQ,R\rangle =\langle P, RQ^T\rangle = \langle Q, P^TR\rangle.
\end{align}
From now on we fix $d=3$.  
Let $\Omega\subset\R^3$ be a bounded domain. Assumptions on the boundary $\partial\Omega$ will be specified where they are needed. For a function $u=(u^1,\ldots, u^m)^T:\Omega\to\R^m$ (with $m\in \N$), the  differential  $\nabla u$ is given by
\begin{align*}
 \nabla u
 =\begin{pmatrix}
 \nabla u^1\\
 \vdots\\
 \nabla u^m
 \end{pmatrix}
 \in \R^{m\times d}\,,
\end{align*}
with $(\nabla u)_{k,l}= \partial_{x_l} u^k$ for $1\leq k\leq m$ and $1\leq \ell\leq d$ and with $\nabla u^k=(u^k_{,x_1},\ldots , u^k_{,x_d})\in \R^{1\times d}$,  while for scalar functions $u:\Omega\to\R$ we also define $\grad u:=(\nabla u)^T\in \R^{d\times 1}$.   
For a vector field $w:\Omega\to\R^3$, the divergence and the  curl are  given as
\begin{align*}
\dyw w=\sum_{i=1}^ 3 w^{i}_{, x_i},\qquad 
\curl w =\big(w^{3}_{,x_2}-w^{2}_{,x_3},w^{1}_{,x_3}-w^{3}_{,x_1},w^{1}_{,x_2}-w^{2}_{,x_1}\big)\,.
\end{align*}
For tensor fields $Q:\Omega\to\R^{n\times 3}$ ($n\in \N$), $\Curl Q$ and $\Dyw Q$ are defined row-wise: 
\begin{align*}
 \Curl Q=\begin{pmatrix} \curl Q^1
             \\
             \vdots\\
             \curl Q^n    
            \end{pmatrix}
\in  \R^{n\times 3}\,,\qquad 
\text{ and } 
 \Dyw Q=\begin{pmatrix}
  \dyw Q^1
  \\ \vdots \\
   \dyw Q^n
   \end{pmatrix}\in \R^n\,,
\end{align*}
where $Q^i$ denotes the $i$-th row of $Q$. With these definitions, for $u:\Omega\to \R^m$ we have consistently $ \Curl\nabla u=0\in \R^{m \times d}$.  

$C_0^{\infty}(\Omega)$ is the set of smooth functions with compact support in $\Omega$. Additionally,  $L^2(\Omega)$ denotes the usual Lebesgue spaces of square integrable scalar functions, vector or tensor fields on $\Omega$ with values in $\R$, $\R^d$ or $\R^{m\times n}$, respectively. If needed, we will write more explicitly $L^2(\Omega,\R^s)$ to indicate that we deal with $\R^s$-valued functions. The Sobolev spaces \cite{adamssobolev} used in this paper are
\begin{align*}
& H^1(\Omega)=\{u\in L^2(\Omega)\, |\,\, \nabla u\in L^2(\Omega)\}\,, \qquad \|u\|^2_{H^1(\Omega)}:=\|u\|^2_{L^2(\Omega)}+\|\nabla u\|^2_{L^2(\Omega)}\, ,\\[2ex]
&H({\curl};\Omega)=\{v\in L^2(\Omega;\R^d)\, |\,\, \curl v\in L^2(\Omega)\}\,,\qquad\|v\|^2_{H({\rm curl};\Omega)}:=\|v\|^2_{L^2(\Omega)}+\|{\rm curl}\, v\|^2_{L^2(\Omega)}\, ,\\[2ex]
&H(\dyw;\Omega)=\{v\in L^2(\Omega;\R^d)\, |\,\, \dyw v\in L^2(\Omega)\},\qquad
\|v\|^2_{H({\dyw};\Omega)}:=\|v\|^2_{L^2(\Omega)}+\|\dyw v\|^2_{L^2(\Omega)}\, .
\end{align*}
Moreover, $H_0^1(\Omega)$ is the completion of $C_0^{\infty}(\Omega)$ with respect to the $H^1$-norm and $H_0({\rm curl};\Omega)$ and $H_0(\dyw;\Omega)$ are the completions of $C_0^{\infty}(\Omega)$ with respect to the $H(\curl)$-norm and the $H(\dyw)$-norm, respectively. 
By $H^{-1}(\Omega)$ we denote the dual of $H_0^1(\Omega)$. 
  Finally, 
\begin{equation*}
    L^2_{{\rm loc}}(\Omega)=\{u:\Omega\rightarrow\R\, |\,\,  u\in L^2(\tilde{\Omega})\,\,\mathrm{for\, all}\,\, \tilde{\Omega}\Subset\Omega\}\,,
\end{equation*}
where $\tilde{\Omega}\Subset\Omega$ means that there exists $K$ compact such that  $\tilde{\Omega}\subset K\subset\Omega$. This definition extends in a natural way to  the above Sobolev spaces.

\section{Local regularity via generalized finite differences}
\label{s:regularity}
In a first step we study local regularity on a cube. 
For $r>0$ and $m,n\in \N$,  $d=3$,  let $\Omega:=C_r:=(-r,r)^d$ denote the
cube with side length $2r$ centered at the origin. 

Given a   coefficient tensor
 $\bbA:\Omega \rightarrow \Lin(\R^{m}\times \R^{n\times d}\times
  \R^{m\times d} \times \R^{n\times d},\R^{m}\times \R^{n\times d}\times
  \R^{m\times d} \times \R^{n\times d}) $, for $u,v\in H^1(\Omega,\R^m)$ and 
$P,Q\in H(\Curl;\Omega,\R^{n\times d})$ we define the bilinear form
\begin{align}
\label{mod_def_bili}
a(( u,P);( v,Q) )
:=\int_\Omega \Big\langle 
\bbA(x)\begin{pmatrix}u\\ P \\ \nabla u\\ \Curl P \end{pmatrix}, 
\begin{pmatrix}v\\ Q \\ \nabla v\\ \Curl Q \end{pmatrix}
\Big\rangle
\dx \,.
\end{align}
It is assumed that the coefficient tensor $\bbA$ and the associated bilinear form satisfy
\begin{itemize}
\item[(A1)] $\bbA \in C^{0,1}(\overline{\Omega})$ 
and there exist constants $c\in \R$ and 
  $\alpha>0$  such that for all $u\in H_0^1(\Omega,\R^m)$ and $P\in
  H_0(\Curl;\Omega,\R^{n\times d})$ it holds
\begin{align}
\label{mod_ineq_coerc}
a(( u;P);  (u,P)) 
\geq \alpha\big(\norm{\nabla u}^2_{L^2(\Omega)} + \norm{P}^2_{L^2(\Omega)} +
\norm{\Curl P}^2_{L^2(\Omega)}\big) - c\norm{u}^2_{L^2(\Omega)}.
\end{align}
\end{itemize}
The right hand sides shall be chosen as follows:
\begin{itemize}
\item[(A2)] $f\in L^2_\text{loc}(\Omega,\R^m)\cap H^{-1}(\Omega,\R^m)$\,, $M\in H^1_\text{loc}(\Omega;\R^{n\times d})\cap \big(H_0(\Curl;\Omega,\R^{n\times d})\big)^{\ast}$, where \\$ \big(H_0(\Curl;\Omega,\R^{n\times d})\big)^{\ast}$ is the dual space of the Hilbert space $H_0(\Curl;\Omega,\R^{n\times d})$.
\end{itemize}
\begin{theorem}[Local regularity on a cube]
\label{them_reg_cube}
Let $\Omega=C_r$ and assume (A1) and (A2). Let $u\in H^1(\Omega,\R^m)$ and 
$P\in H(\Curl; \Omega,\R^{n\times d})$ satisfy
\begin{align}
\label{def_weak_form1}
a((u,P);(v,Q))= \int_\Omega \langle f,v\rangle \dx + \int_\Omega\langle M,Q\rangle\dx
\end{align}
for all $v\in H_0^1(\Omega,\R^m)$ and all $Q\in H_0(\Curl;\Omega,\R^{n\times d})$. 
Then $u\in H^2_\text{loc}(\Omega;\R^m)$ and $P,\Curl P\in H^1_\text{loc}(\Omega;\R^{n\times d})$. 
\end{theorem}

The proof is based on a modified difference quotient argument. As already mentioned in the introduction,  the challenge
is to choose the test functions (finite differences of the solution $(u,P)$ with
a suitable cut-off) 
 in such a way that they belong to the test space 
$ H_0^1(\Omega)\times  H_0(\Curl;\Omega)$ and that the $\curl$-structure is maintained, see \eqref{id_curl_schoeberl_weak}. In particular for the $P$
component we have to modify the classical ansatzes. Instead of standard finite differences we work with differences based on inner variations and additionally transform the tensor valued  function $P$ by a $\Curl$-adapted version of the Piola-transform. 

\subsection{Preliminary information  about inner variations} 

Let $\rho\in (0,r)$ and choose $\varphi\in C_0^\infty(C_r)$ with
$\supp\varphi\subset C_\rho$ and $\varphi(x)=1$ for $x\in C_{\rho/2}$. For
$h\in \R^d$ we define the following family of inner variations
\begin{align*}
T_h:\R^d\rightarrow \R^d\,,\qquad T_h(x) = x+ \varphi(x) h.
\end{align*}
Obviously, we have $T_h(x)=x$ for $x\in \R^d\backslash C_\rho$ and 
\begin{align*}
\nabla T_h(x) = \idd + h\otimes \grad \varphi(x)\,,\qquad 
\det \nabla T_h(x)= 1 + \langle h, \grad \varphi(x)\rangle.
\end{align*}
There exists $h_0>0$ such that for every 
$ h \in \R^d$  with $\abs{h}\leq h_0$,
 the
mapping $T_h$ is a diffeomorphism from $C_r$
 onto itself, see
e.g. \cite{var:GH96a}. 
For every $h\in \R^d$ with $\abs{h}\leq h_0$ and every $x\in \R^d$  it further holds
\begin{equation}
\label{invers}
(\nabla T_h)^{-1}(x)=\idd-\big(1+\langle h,\grad\varphi(x)\rangle\big)^{-1}h\otimes\grad\varphi(x)\,.
\end{equation}
Moreover, the following uniform estimate is valid:
\begin{align}
\label{est:Th}
 \sup_{\abs{h}\leq h_0}\Big(\norm{\det\nabla T_h}_{L^\infty(\R^d)} + 
 \norm{\det(\nabla T_h)^{-1}}_{L^\infty(\R^d)} + \norm{ T_h}_{W^{1,\infty}(\R^d)} + 
 \norm{ T_h^{-1}}_{W^{1,\infty}(\R^d)}
 \Big)=:C_{h_0}<\infty\,.
\end{align}
%
For a function $Q:\Omega\rightarrow \R^{n\times d}$ we define the Piola-type transformation 
\begin{align}
\label{def:Piolatypetrafo}
\calT_h(Q)(x):=  Q(T_h(x)) \nabla T_h(x)
\end{align}
and set $S_h:=T_h^{-1}$ and $\calS_h(Q):=  Q\circ S_h\nabla S_h$ if $\abs{h}\leq h_0$. 
 
\begin{lemma}
\label{lem:piolacurl}
For every $h\in \R^d$ with $\abs{h}\leq h_0$ the mapping
 $\calT_h:
H(\Curl;\Omega,\R^{n\times d})\rightarrow H(\Curl;\Omega,\R^{n\times d})$ is well defined, linear and continuous and there exists a constant $c>0$ such that for all $Q\in H(\Curl;\Omega)$
\begin{align}
\label{est_trafo_schoeberl}
\sup_{\abs{h}\leq h_0} \big(\norm{\calT_h(Q)}_{L^2(\Omega)} + \norm{\Curl
  \calT_h(Q)}_{L^2(\Omega)}\big) \leq 
c\big( 
\norm{Q}_{L^2(\Omega)} + \norm{\Curl
  Q}_{L^2(\Omega)}\big). 
\end{align}
Moreover, for every $Q\in H(\Curl;\Omega,\R^{n\times d})$ we have
\begin{align}
 \label{id_curl_schoeberl_weak}
\Curl\big( \calT_h(Q)\big) = (\det\nabla T_h)  (\Curl Q)\circ T_h \big(\nabla T_h\big)^{-T}\,.
 \end{align}
\end{lemma}
\begin{proof}
Let $h\in \R^d$ with $\abs{h}\leq h_0$. From \eqref{eq:schoeberl_strong}, 
 see for instance \cite{Schoberlnotes}, we know that  for every $\psi\in C_0^\infty(\Omega;\R^{n\times d})$ the identity \eqref{id_curl_schoeberl_weak} is satisfied for both, $\calT_h(\psi)$ as well as $\calS_h(\psi)$.  Hence, for the distributional curl of $Q\in H(\Curl;\Omega)$ and every $\psi\in C_0^\infty(\Omega;\R^{n\times d})$ we find after a transformation of coordinates ($y=T_h(x)$, $x=S_h(y)$) and taking into account \eqref{eq:sp_matrix} and the identity   \eqref{id_curl_schoeberl_weak}
\begin{align*}
 \langle \Curl \calT_h( Q),\psi\rangle&=\int_\Omega \big\langle\calT_h(Q),  \Curl\psi\big\rangle\,\dx \\
 &=\int_\Omega\big\langle Q,(\Curl\Psi)\circ S_h (\nabla S_h)^{-T}\det\nabla S_h\big\rangle\,\dy
 \\
 &= \int_\Omega \big\langle Q,\Curl(\calS_h(\psi))\big\rangle\,\dy
 =\int_\Omega \big\langle \Curl Q,\calS_h(\psi)\big\rangle\,\dy
 \\
 &=\int_\Omega \big\langle \det \nabla T_h  (\Curl Q)\circ T_h\, (\nabla T_h)^{-T},\psi\big\rangle\,\dx\,.
\end{align*}
Since by assumption $\Curl Q$ (and hence also $(\Curl Q)\circ T_h$) belongs to $L^2(\Omega)$,  these calculations show that $\Curl(\calT_h(Q))\in L^2(\Omega)$ as well and that $\eqref{id_curl_schoeberl_weak}$ is satisfied. Estimate \eqref{est_trafo_schoeberl} is an immediate consequence of the identity \eqref{id_curl_schoeberl_weak} combined with  \eqref{est:Th}.
\end{proof}
 With
 $\rho\in (0,r)$, $\varphi\in C_0^\infty(C_r)$ and $T_h,S_h$  defined as above, 
for $\abs{h}\leq h_0$ and $w:\Omega\rightarrow\R^s$ we set
\begin{align*}
\triangle_{T_h}w:= w - w\circ T_h\,,\qquad 
\triangle_{S_h}w:= w - w\circ S_h(x)\,.
\end{align*}
The next lemma recalls estimates for finite differences of functions from $H^1(\Omega)$:
\begin{lemma}
\label{lem:diffquotH1}
 There exists a constant $c>0$ such that for every $h\in \R^d$ with $\abs{h}\leq h_0$ and every $w\in H^1(\Omega;\R^m)$ we have
 \begin{align}
 \label{est:h1diff}
  \norm{\triangle_{T_h}w}_{L^2(\Omega)} + \norm{\triangle_{S_h}w}_{L^2(\Omega)} \leq c\abs{h}\norm{w}_{H^1(C_\rho)}\,.
 \end{align}
Moreover, let $\triangle_{S_h}^* : H^{-1}(\Omega,\R^m)\rightarrow H^{-1}(\Omega,\R^m)$ denote the adjoint   operator of $\triangle_{S_h}$ given by  $\langle\triangle_{S_h}^*(f),v\rangle:= \langle f,
\triangle_{S_h}v\rangle$ for all $v\in H_0^1(\Omega,\R^m)$. Then there exists a constant $c>0$ such that for every $f\in L^2_\text{loc}(\Omega;\R^m)\cap H^{-1}(\Omega,\R^m)$ and every $h\in \R^m$ with $\abs{h}\leq h_0$ it holds
\begin{align}
 \label{est:dualnorm}
 \norm{\triangle_{S_h}^*(f)}_{H^{-1}(\Omega;\R^m)}\leq c\abs{h}\norm{f}_{L^2(C_\rho)}.
\end{align}
\end{lemma}
\begin{proof}
 Estimate \eqref{est:h1diff} follows by density of $C^\infty(\overline{\Omega})$ in $H^1(\Omega)$ and taking into account that $\supp(\text{id} - T_h)\subset \overline{C_\rho}$. Estimate  \eqref{est:dualnorm} is immediate by duality.
\end{proof}

\subsection{Proof of Theorem \ref{them_reg_cube}}
Let  $(u,P)\in H^1(\Omega,\R^m)\times H(\Curl;\Omega,\R^{n\times d})$ be given as in
Theorem \ref{them_reg_cube}. We define 
\begin{align}
u_h&:= \triangle_{S_h}(\triangle_{T_h} u),
\label{def_testu}\\
P_h&:= (1 - \calS_h)((1-\calT_h) P).
\label{def_testq}
\end{align}
Clearly, $u_h\in H^1(\Omega;\R^m)$ and thanks to  Lemma \ref{lem:piolacurl} we have $P_h\in H(\Curl;\Omega,\R^{n\times d})$. Moreover, since $T_h(x)=x$ on $\R^d\backslash C_\rho$,
it holds $\supp(u_h,P_h)\subset \overline{C_\rho}$ which implies that $(u_h,P_h)\in H_0^1(\Omega,\R^m)\times
H_0(\Curl;\Omega,\R^{n\times d})$. Hence, these functions are admissible test functions
for the weak formulation \eqref{def_weak_form1} implying that for all
$\abs{h}\leq h_0$ we have 
\begin{align}
\label{def_weak_form2}
a\big((u,P);(u_h,P_h)\big)=\langle f,u_h\rangle + \langle M,P_h\rangle.
\end{align}
The aim is to derive estimate \eqref{eq:regestdiffquot}, from which the regularity result ensues. 
The first term on the right hand side of  \eqref{def_weak_form2}  
can be  estimated based on assumption (A2) and Lemma \ref{lem:diffquotH1}: 
\begin{align}
\langle f,u_h\rangle_{H^{-1}(\Omega),H^1_0(\Omega)}
&= \langle \triangle_{S_h}^* f, 
\triangle_{T_h}u\rangle_{H^{-1}(\Omega),H^1_0(\Omega)} 
\nonumber\\
&\leq \norm{\triangle_{S_h}^* f}_{H^{-1}(\Omega)}
\norm{\triangle_{T_h}u}_{H^1(\Omega)} 
\nonumber\\
&\leq c \abs{h}\norm{f}_{L^2(C_\rho)}\norm{\triangle_{T_h}u}_{H^1(\Omega)} .
\label{est_f1}
\end{align}
Let us emphasize that the constant $c>0$ does not depend on $h$ or $f$ (but of course on $\rho$,
$\varphi$).  

The second term on the right hand side of  \eqref{def_weak_form2}  
 is estimated in a similar
 way using
again (A2): With $Q_h:=(1-\calT_h)P$ and after a transformation of coordinates
\begin{align*}
\langle M,(1-\calS_h)(1-\calT_h) P
\rangle_{\big(H_0(\Curl;\Omega)\big)^{\ast},H_0(\Curl;\Omega)} 
&= \int_\Omega \big\langle M-(\det\nabla T_h)  (M\circ T_h)(\nabla T_h)^{-T},Q_h\big\rangle\,\dx
\end{align*}
The assumption that $M\in H^1_\text{loc}(\Omega)$ and the structure of $T_h$ imply that 
\begin{align*}
 &\norm{ M-(\det\nabla T_h)  (M\circ T_h)(\nabla T_h)^{-T}}_{L^2(\Omega)}
 \nonumber \\
 &\qquad \leq
 \norm{M - M\circ T_h}_{L^2(C_\rho)} + \norm{M\circ T_h \big(\idd - \det\nabla T_h (\nabla T_h)^{-T}\big)}_{L^2(C_\rho)}
 \nonumber \\
 &\qquad\leq 
 c\abs{h}\norm{M}_{H^1(C_\rho)}
\end{align*}
with a constant $c$ that is independent of $h$ and $M$. 
Hence, 
\begin{align}
\label{est:M}
\abs{\langle M,(1-\calS_h)(1-\calT_h)
P\rangle_{H_0(\Curl;\Omega)^*,H_0(\Curl;\Omega)}} 
\leq c\abs{h} \norm{M}_{H^1(C_\rho)}\norm{(1-\calT_h)P}_{L^2(\Omega)}\,.
\end{align}
In order to estimate the left hand side of 
\eqref{def_weak_form2} let
$w_h:=\triangle_{T_h}u$ and $Q_h:=(1-\calT_h)P$. 
Observe that $u_h=\triangle_{S_h}(\triangle_{T_h}u)=w_h - w_h\circ S_h$ and $ P_h= (1 - \calS_h)((1-\calT_h) P) = Q_h - \calS_h Q_h$. 
Hence,
\begin{align}
\label{id_pr1}
a((u,P), (u_h,P_h))
&= a((u,P);(w_h,Q_h))-a((u,P);(w_h\circ S_h,\calS_h Q_h)).
\end{align}
In the next step, 
 we rewrite
the second term on the right hand side.  
For  $v,w\in H^1(\Omega,\R^m)$, $Q,R\in H(\Curl;\Omega,\R^{n\times d})$, we define 
\begin{align*}
 a_h((v,R);(w,Q))
&=
\int_\Omega
\det\nabla T_h 
\Big\langle
\bbA\circ T_h 
\begin{pmatrix} 
v\\
R\,  \nabla  T_h^{-1}
\\
\nabla v\, \nabla T_h^{-1}
\\
\frac{1}{\det \nabla T_h} (\Curl R) (\nabla T_h)^T
\end{pmatrix}
,
\begin{pmatrix} 
w\\
  Q \,\nabla T_h^{-1} 
\\
\nabla w \,\nabla T_h^{-1}
\\
\frac{1}{\det \nabla T_h}  (\Curl Q) (\nabla T_h)^{T}
\end{pmatrix}
\Big\rangle\,
\dx
\\
&\equiv \int_{\Omega} \Big\langle\bbA_h(x)
\begin{pmatrix} 
v\\
 R 
\\
\nabla v
\\
\Curl R
\end{pmatrix}
,\begin{pmatrix} 
w\\
 Q 
\\
\nabla w
\\
\Curl Q
\end{pmatrix}
\Big\rangle\,\dx\,.
\end{align*}
Direct calculations based on the transformation $x=S_h(y)$ and on the identity \eqref{id_curl_schoeberl_weak} then yield 
\begin{align}
 a\big((u,P);(w_h\circ S_h,\calS_h(Q_h))\big) &= 
 a_h\big((u\circ T_h,\calT_h(P));(w_h,Q_h)\big)\,.
\end{align}
With this, identity \eqref{id_pr1} can be rewritten as 
\begin{align}
 a((u,P), (u_h,P_h))
&= a((u,P);(w_h,Q_h))-a((u,P);(w_h\circ S_h,\calS_h Q_h))
\nonumber\\
&= a((u,P);(w_h,Q_h)) - a_h\big((u\circ T_h,\calT_h(P));(w_h,Q_h)\big)
\nonumber\\
&= a((w_h;Q_h);(w_h,Q_h)) 
\nonumber\\
&\qquad + a((u\circ T_h,\calT_h(P));(w_h,Q_h))- a_h\big((u\circ T_h,\calT_h(P));(w_h,Q_h)\big)\,.
\label{id_bili_trafo}
\end{align}
Due to the coercivity assumption (A1) the first term is estimated 
as follows 
\begin{align}
\label{est_bili_coerc}
a\big((w_h,Q_h);(w_h,Q_h)\big) \geq \alpha\left( 
\norm{\triangle_{T_h}u}_{H^1(\Omega)}^2 +
\norm{(1-\calT_h)(P)}_{H_{\curl}(\Omega)}^2 \right) -
c
\abs{h}^2\norm{u}_{H^1(\Omega)}^2,
\end{align}
where  we have also used that $u\in H^1(\Omega,\R^m)$. The
constant $c$ is independent of $u,P,h$ but depends on $\rho$ and the chosen
cut-off function $\varphi$.
  
In order to estimate the last two terms in \eqref{id_bili_trafo},
observe first that thanks to assumption (A1) and the special structure of $\nabla T_h$  there exists a constant $c>0$ such that for all $\abs{h}\leq h_0$ we
 have 
\begin{align*}
\norm{\bbA_h - \bbA}_{L^\infty(\Omega)} \leq
c\abs{h}\norm{\bbA}_{C^{0,1}(\overline\Omega)}.
\end{align*}
Hence, 
\begin{align}
&\nquad\nquad
\abs{a\big((u\circ T_h, \calT_h(P));(w_h,Q_h)\big)
- a_h\big((u\circ T_h, \calT_h(P));(w_h,Q_h)\big)}
\nonumber\\
&\leq 
c \abs{h} 
\big(\norm{u}_{H^1(\Omega)} + \norm{P}_{H(\Curl,\Omega)}\big)
\big(\norm{\triangle_{T_h}u}_{H^1(\Omega)} +
 \norm{(1-\calT_h)P}_{H(\Curl,\Omega)}\big),
\label{est_bili_2}
\end{align}
where we also applied estimate \eqref{est_trafo_schoeberl}.

Combining estimates \eqref{est_f1}, \eqref{est:M}, \eqref{est_bili_coerc}
 and \eqref{est_bili_2} and applying
Young's inequality on the right hand side
 (in order to absorb the terms $\norm{\triangle_{T_h}u}_{H^1(\Omega)}$ and 
$\norm{(1-\calT_h)(P)}_{H(\Curl;\Omega)}$ 
into the left hand side) we finally obtain:
\begin{multline}
\label{eq:regestdiffquot}
\norm{\triangle_{T_h}u}_{H^1(\Omega)}^2 +
\norm{(1-\calT_h)(P)}_{H(\Curl;(\Omega)}^2
\\
\leq c\abs{h}^2 \Big(\norm{f}_{L^2(C_\rho)} + \norm{M}_{H^1(C_\rho)}+
\norm{u}_{H^1(\Omega)}  +\norm{P}_{H(\Curl;\Omega)}\Big)^2,
\end{multline}
which is the desired local regularity estimate. Indeed, from the previous estimate it follows that 
\begin{multline*}
 \sup_{\substack{h\in \R^d\backslash\{0\},\\ \abs{h}\leq h_0}} 
 \abs{h}^{-1}\Big(\norm{u(\cdot + h) - u(\cdot)}_{H^1(C_{\rho/2})} + 
 \norm{P(\cdot + h) - P(\cdot)}_{H(\Curl;C_{\rho/2})}\Big) 
 \\ 
 \leq 
 c\Big( \norm{f}_{L^2(C_\rho)} + \norm{M}_{H^1(C_\rho)} +
\norm{u}_{H^1(C_\rho)}  +\norm{P}_{H(\Curl;C_\rho)}\Big),
\end{multline*}
implying that $u\in H^1(C_{\rho/2})$ and $P,\Curl P\in H^1(C_{\rho/2})$, see for instance \cite[Lemma 7.24]{GiTr86}. 

The above reasoning can be repeated for arbitrary cubes 
$C_\rho(x_0):=\{x\in \R^d\, |\,\,\norm{x-x_0}_{\infty}<\rho\}$ that are compactly embedded in $\Omega$, which ultimately finishes the proof of Theorem \ref{them_reg_cube}.

\section{Application to the relaxed micromorphic model}
\label{sec:micro}
The relaxed micromorphic model \cite{NGperspective,Ghibastatic,Ghibadyn,Owczghibaneffexist,exprolingalberdi} is a special variant of the classical Mindlin-Eringen micromorphic elasticity model \cite{Mindlin64,Eringen68}. It features the classical three translational degrees of freedom -- the displacement $u:\Omega\rightarrow\R^3$ -- as well as nine additional degrees of freedom --  the non-symmetric micro-distortion $P:\Omega\rightarrow\R^{3\times 3}$. Contrary to the classical micromorphic model, the curvature part does not consider the full gradient $\nabla P\in\R^{3\times 3\times 3}$, but only the matrix $\Curl P$ (hence the name "relaxed").

Let  $\bbC_e\in \Lin(\Sym(3),\Sym(3))$, $\bbC_\text{micro}\in \Lin(\Sym(3),\Sym(3))$ and $\bbL_c\in \Lin(\R^{3\times 3},\R^{3\times 3})$ all be
symmetric and positive definite. The model reads as follows: given an external body force density $f:\Omega\rightarrow \R^3$ and a body moment tensor $M:\Omega\rightarrow \R^{3\times 3}$ find a displacement field $u:\Omega\to \R^3$ and a  micro-distortion tensor $P:\Omega\to\R^{3\times 3}$ satisfying
\begin{align*}
 0&=\Dyw\Big(\bbC_e\sym(\nabla u - P)\Big) + f\quad \text{in }\Omega,\\
 0&= -\Curl\Big(\bbL_c\Curl P\big) + \bbC_e\sym(\nabla u - P) - \bbC_\text{micro}\sym P + M\quad \text{in }\Omega.
\end{align*}
This system needs to be closed by boundary conditions. In this system, the equations of  linear elasticity are coupled with a matrix valued Maxwell-type problem. The curvature constitutive tensor $\Lc$ introduces a characteristic length in the model such that smaller samples are  stiffer. The interest of such type of modelling derives from the need to describe micro-structured materials with a continuum model. 
Here, we only consider the more challenging case with symmetric force-stress  tensor $\sigma=\bbC_e\sym(\nabla u-P)$, i.e.\ zero Cosserat couple moduli $\bbC_c$. 
In the dynamic case the relaxed micromorphic model has shown its unique capability to phenomenologically describe frequency-band gaps for periodic structures \cite{rizzi2021exploring,rizzi2021boundary}. 
Another strong point of the model, unknown for other generalised continuum models, is that the size-independent tensors $\D$ and $\Ha$ can be directly interpreted  \cite{neff2020identification}. We refer to \cite{Ghibadyn,OGNdynamicreg} for the discussion of dynamic versions of this model. 

The bilinear form associated with this model is given by
\begin{multline}
\label{bili_micromorph}
a((u,P);(v,Q)):=\int_\Omega \langle\bbC_e\sym(\nabla u- P),\sym(\nabla v - Q)\rangle +
\langle \bbC_\text{micro}\sym
P,\sym Q \rangle
\\+\langle \bbL_c\Curl P, \Curl Q\rangle\,\dx,
\end{multline}
and is defined for $u,v\in H^1(\Omega,\R^3)$ and $P,Q\in
H(\Curl;\Omega,\R^{3\times 3})$. This bilinear form satisfies condition (A1) for every bounded Lipschitz domain $\Omega\subset \R^3$,  which is a consequence of \cite{NGperspective} in combination with the incompatible Korn's inequality discussed in \cite{neff2012poincare,NeffLP,lewintan_neff_2021,LMN21}.

Let now $\Omega\subset\R^3$ be a bounded domain with a Lipschitz boundary. 
Thanks to the coercivity of the bilinear form \eqref{bili_micromorph}, by the Lax-Milgram Theorem   for every $f\in H^{-1}(\Omega;\R^3)$ and $M\in \big(H_0(\Curl;\Omega,\R^{3\times 3})\big)^{\ast}$ there exist unique functions $u\in H_0^1(\Omega,\R^3)$ and $P\in H_0(\Curl;\Omega,\R^{3\times 3})$ that satisfy
\[
\forall (v,Q)\in H_0^1(\Omega;\R^3)\times H_0(\Curl;\Omega,\R^{3\times 3})\qquad a\big((u,P);(v,Q)\big) = \langle f,u\rangle + \langle M,Q\rangle\,.
\]
In \cite{sky2021hybrid,schroder2021lagrange,sky2021Primal},  finite element formulations for this problem have been investigated. From the latter derives the need to exactly understand with what kind of regularity beyond $ H^1(\Omega)\times  H(\Curl;\Omega)$ one can work. As an immediate consequence of Theorem \ref{them_reg_cube} one obtains 
\begin{theorem}
 Let $f\in L^2_\text{loc}(\Omega,\R^3)$ and $M\in H^1_\text{loc}(\Omega;\R^{3\times 3})$. Then the weak solution pair $(u,P)$ satisfies $u\in H^2_\text{loc}(\Omega)$ and $P$, $\Curl P \in H^1_\text{loc}(\Omega)$.
\end{theorem}
 This allows, in the smooth case, to just use standard $H^1$-elements, avoiding the more demanding implementation of N\'ed\'elec - $H(\Curl;\Omega)$ elements, see  \cite{schroder2021lagrange,sky2021Primal,skyMuench2022} for a discussion of the latter approach. 

\begin{appendix}
\section{Some formulas from vector calculus}
In this section we collect formulas from vector calculus adapted to the notation/definitions in Section \ref{s:notation}. 
 Let $\Omega,\wt \Omega\subset\R^n$ be  bounded domains, $v:\Omega\to\R$, $f:\Omega\to \R^m$, $A:\Omega\to \R^{ n\times m}$ smooth functions. Then the following product rules involving the divergence operator are valid: 
 \begin{align}
 \label{eq:prod-div}
  \dyw(v f)&= \langle\grad v, f\rangle + v\dyw f, \qquad 
  \dyw(A f)= \langle \Dyw (A^T), f\rangle\, + \langle A^T, \nabla f\rangle.
 \end{align}
For diffeomorphisms $\Phi:\Omega\to \wt\Omega$ and tensor valued functions $A:\Omega\to\R^{k\times n}$ the (contravariant) Piola transform is defined as, \cite{Cia88}: 
\[
 \calP_\Phi(A):\wt \Omega\to\R^{k\times n},\qquad 
 \calP_\Phi(A):= (\det \nabla \Phi) \, A\circ\Phi \, (\nabla \Phi)^{-T}\,,
\]
and we have 
\[
 \Dyw \calP_\Phi(A)= \det \nabla \Phi \,(\Dyw A)\circ\Phi\,,
\]
which is a consequence of the Piola identity $\Dyw \big((\det\nabla \Phi)(\nabla \Phi)^{-T}\big)=0$ and \eqref{eq:prod-div}, \cite{Cia88}. The covariant Piola transform is defined as follows for tensor valued functions $A:\Omega\to\R^{k\times n}$:
\begin{align}
 \calQ_\Phi(A):\wt \Omega\to\R^{k\times n},\qquad 
 \calQ_\Phi(A):= (A\circ\Phi) \nabla \Phi
\end{align}
For $n=3$ (and $k\in \N$) the following identity is valid, \cite{Schoberlnotes}:
\begin{align}
\label{eq:schoeberl_strong}
 \Curl(\calQ_\Phi(A)) = \det \nabla \Phi \,\big((\Curl A)\circ \Phi\big)\, (\nabla \Phi)^{-T}\,
\end{align}
 
\end{appendix}

\subsubsection*{Acknowledgment}
This work was conceived during a visit of Sebastian Owczarek at the Faculty of Mathematics, Universit\"at Duisburg-Essen, Campus Essen at the invitation of  Patrizio Neff. This visit was sponsored by the Grant of the Polish National Science Center MINIATURA 5, 2021/05/X/ST1/00215. Patrizio Neff and Dorothee Knees acknowledge  support  in the framework of the Priority Programme SPP 2256 "Variational Methods for Predicting Complex Phenomena in Engineering Structures and Materials" funded by the Deutsche Forschungsgemeinschaft (DFG, German research foundation): P.\ Neff within the project "A variational scale-dependent transition scheme - from Cauchy elasticity to the relaxed micromorphic continuum" (Project-ID 440935806), D.\ Knees within the project  "Rate-independent systems in solid mechanics: physical properties, mathematical analysis, efficient numerical algorithms"  (Project-ID 441222077).

\begin{footnotesize}

\def\cprime{$'$} \def\cprime{$'$}

\end{footnotesize}

\end{document}